\documentstyle[12pt]{article}
\newcommand{\be}{\begin{equation}}
\newcommand{\ee}{\end{equation}}
\newcommand{\bea}{\begin{eqnarray}}
\newcommand{\eea}{\end{eqnarray}}
\newcommand{\ba}{\begin{array}}
\newcommand{\ea}{\end{array}}

\newcommand{\bc}{\begin{center}}
\newcommand{\ec}{\end{center}}
\newcommand{\ben}{\begin{enumerate}}
\newcommand{\een}{\end{enumerate}}
\newcommand{\bfi}{\begin{figure}}
\newcommand{\efi}{\end{figure}}

\newcommand{\bq}{\begin{quote}}
\newcommand{\eq}{\end{quote}}
\newcommand{\bqu}{\begin{quotation}}
\newcommand{\equ}{\end{quotation}}
\newenvironment{emphit}{\begin{itemize}}{\end{itemize}}
\newcommand{\bemp}{\begin{emphit}}
\newcommand{\eemp}{\end{emphit}}

\newcommand{\bt}{\begin{tabular}}
\newcommand{\et}{\end{tabular}}

\newtheorem{myth}{Theorem}[section]
\newtheorem{mylem}{Lemma}[section]
\newtheorem{mycor}{Corollary}[section]

\newtheorem{mydef}{Definition}[section]
\newtheorem{myrem}{Remark}[section]

\begin{document}
\date{}
\title{On central loops and the central square property
}
\author{T\`em\'it\d {\'o}p\d {\'e} Gb\d {\'o}l\'ah\`an Jaiy\'e\d
ol\'a$^1$\\
 \&\\
John Ol\'us\d ol\'a Ad\'en\'iran$^2$\thanks{Corresponding author}}
\maketitle
\begin{abstract}
The representation sets of a central square C-loop are investigated.
Isotopes of central square C-loops of exponent $4$ are shown to be
both C-loops and A-loops.
\end{abstract}
\section{Introduction}

C-loops are one of the least studied loops. Few publications that
have considered C-loops include Fenyves \cite{phd50}, \cite{phd56},
Beg \cite{phd169}, \cite{phd170}, Phillips et. al. \cite{phd9},
\cite{phd58}, \cite{phd59}, \cite{phd22}, Chein \cite{phd54} and
Solarin et. al. \cite{phd55}, \cite{phd52}, \cite{phd53},
\cite{phd10}. The difficulty in studying them is as a result of the
nature of their identities when compared with other Bol-Moufang
identities(the element occurring twice on both sides has no other
element separating it from itself). Latest publications on the study
of C-loops which has attracted fresh interest on the structure
include \cite{phd9}, \cite{phd58}, and \cite{phd59}.

LC-loops, RC-loops and C-loops are loops that satisfies the
identities  $(xx)(yz)=(x(xy))z$ , (zy)$(xx)=z((yx)x)$ and
$x(y(yz))=((xy)y)z$ respectively. Fenyves' work in \cite{phd56} was
completed in \cite{phd9}. Fenyves proved that LC-loops and RC-loops
are defined by three equivalent identities. But in \cite{phd9} and
\cite{phd61}, it was shown that LC-loops and RC-loops are defined by
four equivalent identities. Solarin \cite{phd53} named the fourth
identities left middle(LM-) and right middle(RM-) identities and
loops that obey them are called LM-loops and RM-loops respectively.
These terminologies were also used in \cite{phd51}. Their basic
properties are found in \cite{phd58}, \cite{phd56} and \cite{phd49}.
\begin{mydef}
A set $\Pi $ of permutations on a set $L$ is the representation of a
loop $(L,\cdot )$ if and only if
\begin{description}
\item[(i)] $I\in \Pi $ (identity mapping),
\item[(ii)] $\Pi $ is transitive on $L$(i.e for all $x,y\in L$, there exists a unique $\pi \in \Pi$ such that $x\pi =y$),
\item[(iii)] if $\alpha ,\beta \in \Pi$ and $\alpha \beta^{-1}$ fixes one element
of $L$, then $\alpha =\beta $.
\end{description}
The left(right) representation of a loop $L$ is denoted by
$\Pi_\lambda (L)\Big(\Pi_\rho (L)\Big)$ or $\Pi_\lambda(\Pi_\rho)$
and is defined as the set of all left(right) translation maps on the
loop i.e if $L$ is a loop, then
\begin{displaymath}
\Pi_\lambda=\{L_x:L\to L~\vert~ x\in L\}~\textrm{and}~
\Pi_\rho=\{R_x:L\to L~\vert~x\in L\}~\textrm{where}~R_x:L\to
L~\textrm{and}\end{displaymath}
\begin{displaymath}
L_x:L\to L~\textrm{defined
as}~yR_x=yx~\textrm{and}~yL_x=xy~\textrm{respectively for all}
~x,y\in L~\textrm{are bijections}.
\end{displaymath}
\end{mydef}

\begin{mydef}
Let $(L, \cdot )$ be a loop. The left nucleus of $L$ is the set
\begin{displaymath}
N_\lambda (L, \cdot )=\{a\in L : ax\cdot y=a\cdot xy~\forall~x, y\in
L\}.
\end{displaymath}
The right nucleus of $L$ is the set
\begin{displaymath}
N_\rho (L, \cdot )=\{a\in L : y\cdot xa=yx\cdot a~\forall~ x, y\in
L\}.
\end{displaymath}
The middle nucleus of $L$ is the set
\begin{displaymath}
N_\mu (L, \cdot )=\{a\in L : ya\cdot x=y\cdot ax~\forall~x, y\in
L\}.
\end{displaymath}
The nucleus of $L$ is the set
\begin{displaymath}
N(L, \cdot )=N_\lambda (L, \cdot )\cap N_\rho (L, \cdot )\cap N_\mu
(L, \cdot ).
\end{displaymath}
The centrum of $L$ is the set
\begin{displaymath}
C(L, \cdot )=\{a\in L : ax=xa~\forall~x\in L\}.
\end{displaymath}
The center of $L$ is the set
\begin{displaymath}
Z(L, \cdot )=N(L, \cdot )\cap C(L, \cdot ).
\end{displaymath}
\end{mydef}
$L$ is said to be a centrum square loop if~ $x^2\in C(L, \cdot )$
for all $x\in L$. $L$ is said to be a central square loop if~
$x^2\in Z(L, \cdot )$ for all $x\in L$. $L$ is said to be left
alternative if for all $x, y\in L,~ x\cdot xy=x^2y$ and is said to
right alternative if for all $x, y\in L,~ yx\cdot x=yx^2$. Thus, $L$
is said to be alternative if it is both left and right alternative.
The triple $(U, V, W)$ such that $U, V, W\in SYM(L, \cdot )$ is
called an autotopism of $L$ if and only if
\begin{displaymath}
xU\cdot yV=(x\cdot y)W~\forall ~x, y\in L.
\end{displaymath}
$SYM(L, \cdot )$ is called the permutation group of the loop
$(L,\cdot )$. The group of autotopisms of $L$ is denoted by $AUT(L,
\cdot )$. Let $(L, \cdot )$ and $(G, \circ )$ be two distinct loops.
The triple $(U, V, W) : (L, \cdot )\to (G, \circ )$ such that $U, V,
W : L\to G$ are bijections is called a loop isotopism if and only if
\begin{displaymath}
xU\circ yV=(x\cdot y)W~\forall ~x, y\in L.
\end{displaymath}
In \cite{phd60}, the three identities stated in \cite{phd56} were
used to study finite central loops and the isotopes of central
loops. It was shown that in a finite RC(LC)-loop $L$,
$\alpha\beta^2\in \Pi_\rho(L)\big(\Pi_\lambda(L)\big)$ for all
$\alpha,\beta\in \Pi_\rho(L)\big(\Pi_\lambda(L)\big)$ while in a
C-loop $L$, $\alpha^2\beta\in \Pi_\rho(L)\big(\Pi_\lambda(L)\big)$
for all $\alpha,\beta\in \Pi_\rho(L)\big(\Pi_\lambda(L)\big)$. A
C-loop is both an LC-loop and an RC-loop (\cite{phd56}), hence it
satisfies the formal. Here, it will be shown that LC-loops and
RC-loops satisfy the later formula.

Also in \cite{phd60}, under a triple of the form
$(A,B,B)\big((A,B,A)\big)$, alternative centrum square loop
isotopes of centrum square C-loops were shown to be C-loops. It
will be shown here that the same result is true for RC(LC)-loops.

It is shown that a finite loop is a central square central loop if
and only if its left and right representations are closed relative
to some left and right translations.

Central square C-loops of exponent 4 are shown to be groups, hence
their isotopes are both C-loops and A-loops.

For definition of concepts in theory of loops readers may consult
\cite{phd41}, \cite{phd51} and \cite{phd3}.

\section{Preliminaries}

\begin{mydef}\label{definition:bijection}(\cite{phd3})
Let $(L, \cdot )$ be a loop and $U, V, W\in SYM(L, \cdot )$.
\begin{enumerate}
\item If $(U, V, W)\in AUT(L, \cdot )$ for some $U, V, W$, then $U$
is called an autotopic,
\begin{itemize}
\item the set of autotopic bijections in a loop $(L,\cdot )$ is
represented by $\Sigma (L,\cdot )$.
\end{itemize}
\item If there exists $V\in SYM(L, \cdot )$ such that $xU\cdot y=x\cdot
yV$ for all $x, y\in L$, then $U$ is called $\mu $-regular while
$U'=V$ is called its adjoint.
\begin{itemize}
\item The set of all $\mu $-regular bijections in a loop $(L,
\cdot )$ is denoted by $\Phi (L, \cdot )$, while the collection of
all adjoints in the loop $(L, \cdot )$ is denoted by $\Phi ^*(L,
\cdot )$.
\end{itemize}
\end{enumerate}
\end{mydef}

\begin{myth}(\cite{phd3})
If two quasigroups are isotopic then their groups of autotopisms are
isomorphic.
\end{myth}

\begin{myth}\label{lambdarhophi:subgroup}(\cite{phd3})
The set $\Phi (Q,\cdot )$ of all $\mu $-regular bijections of a
quasigroup $(Q,\cdot )$ is a subgroup of the group $\Sigma
(Q,\cdot )$ of all autotopic bijections of $(Q,\cdot )$.
\end{myth}

\begin{mycor}\label{lambdarhophi:isomorphism}(\cite{phd3})
If two quasigroups $Q$ and $Q'$ are isotopic, then the corresponding
groups $\Phi $ and $\Phi '[\Phi ^*$ and $\Phi '^*]$ are isomorphic.
\end{mycor}

\begin{mydef}\label{0:0.6}
A loop $(L,\cdot)$ is called a left inverse property loop or right
inverse property loop (L.I.P.L. or R.I.P.L.) if and only if it obeys
the left inverse property or right inverse property(L.I.P or R.I.P):
$ x^\lambda (xy) = y ~\textrm{or}~(yx) x^\rho=y$. Hence, it is
called an inverse property loop (I.P.L.) if and only if it has the
inverse property (I.P.) i.e. it is both an L.I.P. and an R.I.P.
loop.
\end{mydef}

Most of our results and proofs, are stated and written in dual
form relative to RC-loops and LC-loops. That is, a statement like
'LC(RC)-loop... A(B)' where 'A' and 'B' are some equations or
expressions simply means 'A' is for LC-loops while 'B' is for
RC-loops. This is done so that results on LC-loops and RC-loops
can be combined to derive those on C-loops. For instance an
LC(RC)-loop is a L.I.P.L.(R.I.P.L) loop while a C-loop in an
I.P.L. loop.

\section{Finite Central Loops}
\begin{mylem}\label{0:1}
Let $L$ be a loop. $L$ is an LC(RC)-loop if and only if $\beta\in
\Pi_\rho(\Pi_\lambda)$ implies $\alpha\beta\in\Pi_\rho(\Pi_\lambda)$
for some $\alpha \in \Pi_\rho(\Pi_\lambda)$.
\end{mylem}
{\bf Proof}\\
$L$ is an LC-loop if and only if $x\cdot(y\cdot yz)=(x\cdot yy)z$
for all $x,y,z\in L$ while $L$ is an RC-loop if and only if
$(zy\cdot y)x=z(yy\cdot x)$ for all $x,y,z\in L$. Thus, $L$ is an
LC-loop if and only if $xR_{y\cdot yz}=xR_{y^2}R_z$ if and only if
$R_{y^2}R_z=R_{y\cdot yz}$ for all $y,z\in L$ and $L$ is an RC-loop
if and only if $xL_{zy\cdot y}=xL_{y^2}L_z$ if and only if
$L_{zy\cdot y}=L_{y^2}L_z$. With $\alpha =R_{y^2}(L_{y^2})$ and
$\beta=R_z(L_z)$, $\alpha\beta\in\Pi_\rho(\Pi_\lambda)$. The
converse is achieved by reversing the process.$\spadesuit$

\begin{mylem}\label{0:2}
A loop $L$ is an LC(RC)-loop if and only if $\alpha
^2\beta=\beta\alpha^2$ for all $\alpha\in \Pi_\lambda(\Pi_\rho)$ and
$\beta\in \Pi_\rho(\Pi_\lambda)$.
\end{mylem}
{\bf Proof}\\
$L$ is an LC-loop if and only if $x(x\cdot yz)=(x\cdot xy)z$ while
$L$ is an RC-loop if and only if $(zy\cdot x)x=z(yx\cdot x)$. Thus,
when $L$ is an LC-loop, $yR_zL_x^2=yL_x^2R_z$ if and only if
$R_zL_x^2=L_x^2R_z$ while when $L$ is an RC-loop,
$yL_zR_x^2=yR_x^2L_z$ if and only if $L_zR_x^2=R_x^2L_z$. Thus,
replacing $L_x$($R_x$) and $R_z$($L_z$) respectively with $\alpha$
and $\beta$, the result follows. The converse is achieved by doing
the reverse.$\spadesuit$

\begin{myth}\label{0:3}
Let $L$ be a loop. $L$ is an LC(RC)-loop if and only if $\alpha
,\beta\in \Pi_\lambda (\Pi_\rho )$ implies $\alpha^2\beta\in
\Pi_\lambda (\Pi_\rho )$.
\end{myth}
{\bf Proof}\\
$L$ is an LC-loop if and only if $x\cdot(y\cdot yz)=(x\cdot yy)z$
for all $x,y,z\in L$ while $L$ is an RC-loop if and only if
$(zy\cdot y)x=z(yy\cdot x)$ for all $x,y,z\in L$. Thus when $L$ is
an LC-loop, $zL_{x\cdot yy}=zL_y^2L_x$ if and only if
$L_y^2L_x=L_{x\cdot yy}$ while when $L$ is an RC-loop,
$zR_y^2R_x=zR_{yy\cdot x}$ if and only if $R_y^2R_x=R_{yy\cdot x}$.
Replacing $L_y$($R_y$) and $L_x$($R_x$) with $\alpha$ and $\beta$
respectively, we have $\alpha^2\beta\in \Pi_\lambda (\Pi_\rho )$
when $L$ is an LC(RC)-loop. The converse follows by reversing the
procedure.$\spadesuit$

\begin{myth}\label{0:4}
Let $L$ be an LC(RC)-loop. $L$ is centrum square if and only if
$\alpha \in \Pi_\rho (\Pi_\lambda )$ implies $\alpha\beta\in
\Pi_\rho (\Pi_\lambda )$ for some $\beta \in \Pi_\rho (\Pi_\lambda
)$.
\end{myth}
{\bf Proof}\\
By Lemma~\ref{0:1}, $R_{y^2}R_z=R_{y\cdot yz}(L_{y^2}L_z=L_{zy\cdot
y})$. Using Lemma~\ref{0:2}, if $L$ is centrum square,
$R_{y^2}=L_{y^2}$($L_y^2=R_{y^2}$). So: when $L$ is an LC-loop,
$R_{y^2}R_z=L_y^2R_z=R_zL_y^2=R_zR_{y^2}=R_{y\cdot yz}$ while when
$L$ is an RC-loop,
$L_{y^2}L_z=R_y^2L_z=L_zR_{y^2}=L_zL_{y^2}=L_{zy\cdot y}$. Let
$\alpha =R_z(L_z)$ and $\beta =R_{y^2}(L_{y^2})$, then
$\alpha\beta\in \Pi_\rho (\Pi_\lambda )$ for some $\beta \in
\Pi_\rho (\Pi_\lambda )$.

Conversely, if $\alpha\beta\in \Pi_\rho (\Pi_\lambda )$ for some
$\beta \in \Pi_\rho (\Pi_\lambda )$ such that $\alpha =R_z(L_z)$ and
$\beta =R_{y^2}(L_{y^2})$ then $R_zR_{y^2}=R_{y\cdot
yz}$($L_zL_{y^2}=L_{zy\cdot y})$. By Lemma~\ref{0:1},
$R_{y^2}R_z=R_{y\cdot yz}$($L_{zy\cdot y}=L_{y^2}L_z$), thus
$R_zR_{y^2}=R_{y^2}R_z$($L_zL_{y^2}=L_{y^2}L_z$) if and only if
$xz\cdot y^2=xy^2\cdot z$($y^2\cdot zx=z\cdot y^2x$). Let $x=e$,
then $zy^2=y^2z$($y^2z=zy^2$) implies $L$ is centrum square.
$\spadesuit$

\begin{mycor}\label{0:5}
Let $L$ be a loop. $L$ is a centrum square LC(RC)-loop if and only
if
\begin{enumerate}
\item $\alpha\beta\in \Pi_\rho (\Pi_\lambda )$ for all $\alpha\in
\Pi_\rho (\Pi_\lambda )$ and for some $\beta \in \Pi_\rho
(\Pi_\lambda )$, \item $\alpha\beta\in\Pi_\rho(\Pi_\lambda)$ for all
$\beta\in \Pi_\rho (\Pi_\lambda )$ and for some $\alpha \in
\Pi_\rho(\Pi_\lambda)$.
\end{enumerate}
\end{mycor}
{\bf Proof}\\
This follows from Lemma~\ref{0:1} and
Theorem~\ref{0:4}.$\spadesuit$

\section{Isotope of Central Loops}
It must be mentioned that central loops are not conjugacy closed
loops(CC-loops) as concluded in \cite{phd52} or else a study of the
isotopic invariance of C-loops will be trivial. This is because if
C-loops are CC-loops, then a commutative C-loop would be a group
since commutative CC-loops are groups. But from the constructions in
\cite{phd58}, there are commutative C-loops that are not groups. The
conclusion in \cite{phd52} is based on the fact that the authors
considered a loop of units in a central Algebra. This has also been
observed in \cite{phd64}.

\begin{myth}\label{1:1}
Let $(L,\cdot )$ be a loop.
$L$ is an LC(RC)-loop if and only if
$(R_{y^2},L_y^{-2},I)\big((R_y^2,L_{y^2}^{-1},I)\big)\in
AUT(L,\cdot)$ for all $y\in L$.
\end{myth}

{\bf Proof}\\
According to \cite{phd58}, $L$ is an LC-loop if and only if
$x\cdot(y\cdot yz)=(x\cdot yy)z$ for all $x,y,z\in L$ while $L$ is
an RC-loop if and only if $(zy\cdot y)x=z(yy\cdot x)$ for all
$x,y,z\in L$. $x\cdot(y\cdot yz)=(x\cdot yy)z$ if and only if
$x\cdot zL_y^2=xR_{y^2}\cdot z$ if and only if
$(R_{y^2},L_y^{-2},I)\in AUT(L,\cdot )$ for all $y\in L$ while
$(zy\cdot y)x=z(yy\cdot x)$ if and only if $zR^2\cdot x=z\cdot
xL_{y^2}$ if and only if $(R_y^2,L_{y^2}^{-1},I)\big)\in
AUT(L,\cdot)$ for all $y\in L$.$\spadesuit$

\begin{mycor}\label{1:2}
Let $(L,\cdot)$ be an LC(RC)-loop,
$(R_{y^2}L_x^2,L_y^{-2},L_x^2)\big((R_y^2,L_{y^2}^{-1}R_x^2,R_x^2)\big)\in
AUT(L,\cdot)$ for all $x,y\in L$.
\end{mycor}
{\bf Proof}\\
In an LC-loop $L$, $(L_x^2,I,L_x^2)\in AUT(L,\cdot)$ while in an
RC-loop $L$, $(I,R_x^2,R_x^2)\in AUT(L,\cdot)$. Thus by
Theorem~\ref{1:1} : for an LC-loop,
$(R_{y^2},L_y^{-2},I)(L_x^2,I,L_x^2)=(R_{y^2}L_x^2,L_y^{-2},L_x^2)\in
AUT(L,\cdot)$ and for an RC-loop, $(R_y^2,L_{y^2}^{-1},I)(I,
R_x^2,R_x^2)=(R_y^2,L_{y^2}^{-1}R_x^2,R_x^2)\in
AUT(L,\cdot)$.$\spadesuit$

\begin{myth}\label{1:3}
Let $(L,\cdot)$ be a loop. $L$ is a C-loop if and only if $L$ is a
right (left) alternative LC(RC)-loop.
\end{myth}
{\bf Proof}\\
If $(L,\cdot)$ is an LC(RC)-loop, then by Theorem~\ref{1:1},
$(R_{y^2},L_y^{-2},I)\big((R_y^2,L_{y^2}^{-1},I)\big)\in
AUT(L,\cdot)$ for all $y\in L$. If $L$ has the right(left)
alternative property, $(R_y^2,L_y^{-2},I)\in AUT(L,\cdot )$ for all
$y\in L$ if and only if $L$ is a C-loop.$\spadesuit$

\begin{mylem}\label{1:4}
Let $(L,\cdot)$ be a loop. $L$ is an LC(RC, C)-loop if and only if
$R_{y^2}$($R_y^2$,~$R_y^2$)$\in \Phi(L)$ and
$(R_{y^2})^*=L_y^2\Big((R_y^2)^*=L_{y^2},~(R_y^2)^*=L_y^2\Big)\in
\Phi^*(L)$ for all $y\in L$.
\end{mylem}
{\bf Proof}\\
This can be interpreted from Theorem~\ref{1:1}.$\spadesuit$

\begin{myth}\label{1:5}
Let $(G, \cdot )$ and $(H, \circ )$ be two distinct loops. If $G$
is a central square LC(RC)-loop, $H$ an alternative central square
loop and the triple $\alpha =(A, B, B)~\Big(\alpha =(A, B,
A)\Big)$ is an isotopism of $G$ upon $H$, then $H$ is a C-loop.
\end{myth}

{\bf Proof}\\
$G$ is a LC(RC)-loop if and only if $R_{y^2}$($R_y^2$)$\in\Phi (G)$
and $(R_{y^2})^*=L_y^2\Big((R_y^2)^*=L_{y^2}\Big)\in\Phi ^*(G)$ for
all $x\in G$. Using the idea in \cite{phd62} : $L_{xA}'=B^{-1}L_xB$
and $R_{xB}'=A^{-1}R_xA$ for all $x\in G$. Using
Corollary~\ref{lambdarhophi:isomorphism}, for the case of $G$ been
an LC-loop : let $h~:~\Phi(G)\rightarrow \Phi(H)$ and
$h^*~:~\Phi^*(G)\rightarrow \Phi^*(H)$ be defined as
$h(U)=B^{-1}UB~\forall~U\in \Phi(G)$ and
$h^*(V)=B^{-1}VB~\forall~V\in \Phi^*(G)$. This mappings are
isomorphisms. Using the hypothesis,
$h(R_{y^2})=h(L_{y^2})=h(L_y^2)=B^{-1}L_y^2B=B^{-1}L_yBB^{-1}L_yB=L_{yA}'L_{yA}'
=L_{yA}'^2=L_{(yA)^2}'=R_{(yA)^2}'=R_{(yA)}'^2\in \Phi(H)$.
$h^*[(R_{y^2})^*]=h^*(L_y^2)=B^{-1}L_y^2B=B^{-1}L_yL_yB=B^{-1}L_yBB^{-1}L_yB
=L_{yA}'L_{yA}'=L_{yA}'^2\in \Phi^*(H)$. So, $R_y'^2 \in\Phi (H)$
and $(R_y'^2)^*=L_y'^2\in\Phi ^*(H)$ for all $y\in H$ if and only if
$H$ is a C-loop.

For the case of an RC-loop $G$, using $h$ and $h^*$ as above but now
defined as : $h(U)=A^{-1}~UA~\forall~U\in \Phi(G)$ and
$h^*(V)=A^{-1}VA~\forall~V\in \Phi^*(G)$. This mappings are still
isomorphisms. Using the hypotheses,
$h(R_y^2)=A^{-1}R_y^2A=A^{-1}R_yAA^{-1}R_yA=R_{yB}'R_{yB}'=R_{yB}'^2\in
\Phi(H)$.
$h^*[(R_y^2)^*]=h^*(L_{y^2})=h^*(R_{y^2})=A^{-1}R_y^2A=A^{-1}R_yR_yB=B^{-1}R_yBB^{-1}R_yB
=R_{yA}'R_{yA}'=R_{yA}'^2=R_{(yA)^2}'=L_{(yA)^2}'=L_{yA}'^2\in
\Phi^*(H)$. So, $R_y'^2 \in\Phi (H)$ and $(R_y'^2)^*=L_y'^2\in\Phi
^*(H)$ if and only if $H$ is a C-loop.$\spadesuit$

\begin{mycor}\label{1:6}
Let $(G, \cdot )$ and $(H, \circ )$ be two distinct loops. If~ $G$
is a central square left (right) RC(LC)-loop, $H$ an alternative
central square loop and the triple $\alpha =(A, B, B)~\Big(\alpha
=(A, B, A)\Big)$ is an isotopism of $G$ upon $H$, then $H$ is a
C-loop.
\end{mycor}
{\bf Proof} \\
By Theorem~\ref{1:3}, $G$ is a C-loop in each case. The rest of
the proof follows by Theorem~\ref{1:5}.$\spadesuit$

\begin{myrem}\label{1:7}
Corollary~\ref{1:6} is exactly what was proved in \cite{phd60}.
\end{myrem}

\section{Central square C-loops of exponent 4}
For a loop $(L,\cdot)$, the bijection $J~:~L\to L$ is defined by
$xJ=x^{-1}$ for all $x\in L$.
\begin{myth}\label{c:exponent4}
In a C-loop $(L,\cdot )$, if any of the following is true for all
$z\in L$:
\begin{enumerate}
\item $(I,L_z^2,JL_z^2J)\in AUT(L),$ \item $(R_z^2,I,JR_z^2J)\in
AUT(L),$
\end{enumerate}
then, $L$ is a loop of exponent $4$.
\end{myth}
{\bf Proof} \\
\begin{enumerate}
\item If $(I,L_z^2,JL_z^2J)\in AUT(L)$ for all $z\in L$, then : $x\cdot
yL_z^2=(xy)JL_z^2J$ for all $x,y,z\in L$ implies $x\cdot
z^2y=xy\cdot z^{-2}$ implies $z^2y\cdot z^2=y$. Then $y^4=e$. Hence
$L$ is a C-loop of exponent $4$.
\item If $(R_z^2,I,JR_z^2J)\in AUT(L)$ for all $z\in L$, then :
$xR_z^2\cdot y=(xy)JR_z^2J$ for all $x,y,z\in L$ implies
$(xz^2)\cdot y=[(xy)^{-1}z^2]^{-1}$ implies $(xz^2)\cdot
y=z^{-2}(xy)$ implies $(xz^2)\cdot y=z^{-2}x\cdot y$ implies
$xz^2=z^{-2}x$ implies $z^4=e$. Hence $L$ is a C-loop of exponent
$4$.
\end{enumerate}
$\spadesuit$

\begin{myth}\label{c:centralsquare}
In a C-loop $L$, if the following are true for all $z\in L$ :
\begin{enumerate}
\item $(I,L_z^2,JL_z^2J)\in AUT(L),$ \item $(R_z^2,I,JR_z^2J)\in
AUT(L),$
\end{enumerate}
then, $L$ is a central square C-loop of exponent 4.
\end{myth}
{\bf Proof} \\
By the first hypothesis, If $(I,L_z^2,JL_z^2J)\in AUT(L)$ for all
$z\in L$, then : $x\cdot yL_z^2=(xy)JL_z^2J$ for all $x,y,z\in L$
implies $x\cdot z^2y=xy\cdot z^{-2}$.

By the second hypothesis, If $(R_z^2,I,JR_z^2J)\in AUT(L)$ for all
$z\in L$,then : $xR_z^2\cdot y=(xy)JR_z^2J$ for all $x,y,z\in L$
implies $xz^2\cdot y=z^{-2}(xy)$.

Using the two results above and keeping in mind that $L$ is a C-loop
we have :

$x\cdot z^2y=xz^2\cdot y$ if and only if $xy\cdot z^{-2}=z^{-2}\cdot
xy$. Let $t=xy$ then $tz^{-2}=z^{-2}t$ if and only if
$z^2t^{-1}=t^{-1}z^2$. Let $s=t^{-1}$ then $z^2\in C(L,\cdot )$ for
all $z\in L$.

Since $s$ is arbitrary in $L$, then the last result shows that $L$
is centrum square. Furthermore, C-loops have been found to be
nuclear square in \cite{phd58}, thus $z^2\in Z(L,\cdot )$. Hence
$L$ is a central square C-loop. Finally, by
Theorem~\ref{c:exponent4}, $x^4=e$.$\spadesuit$

\begin{myrem}
In \cite{phd58}, C-loops of exponent 2 were found. But in this
section we have further checked for the existence of C-loops of
exponent 4(Theorem~\ref{c:exponent4}). Also, in \cite{phd58} and
\cite{phd56}, C-loops are proved to be naturally nuclear square.
Theorem~\ref{c:centralsquare} gives some conditions under which a
C-loop can be naturally central square.
\end{myrem}

\begin{myth}\label{1:8}
If $A=(U,V,W)\in AUT(L,\cdot )$ for a C-loop $(L,\cdot )$, then
$A_\rho =(V,U,JWJ)\not \in AUT(L,\cdot )$, but $A_\mu
=(W,JVJ,U),A_\lambda =(JUJ,W,V)\in AUT(L,\cdot )$.
\end{myth}
{\bf Proof}\\
The fact that $A_\mu ,A_\lambda \in AUT(L,\cdot )$ has been shown in
\cite{phd41} and \cite{phd3} for an I. P. L. $L$. Let $L$ be a
C-loop. Since C-loops are inverse property loops, $A_\mu
=(W,JVJ,U),A_\lambda =(JUJ,W,V)\in AUT(L,\cdot )$. A C-loop is both
an RC-loop and an LC-loop. So, $(I,R_x^2,R_x^2),(L_x^2,I,L_x^2)\in
AUT(L,\cdot )$ for all $x\in L$. Thus, if $A_\rho\in AUT(L,\cdot )$
when $A=(I,R_x^2,R_x^2)$ and $A=(L_x^2,I,L_x^2)$, $A_\rho =
(I,L_x^2,JL_x^2J)\in AUT(L)$ and $A_\rho =(R_x^2,I,JR_x^2J)\in
AUT(L)$ hence by Theorem~\ref{c:exponent4} and
Theorem~\ref{c:centralsquare}, all C-loops are central square and of
exponent 4(in fact it will soon be seen in Theorem~\ref{1:10} that
central square C-loops of exponent 4 are groups), which is false.
So, $A_\rho =(V,U,JWJ)\not \in AUT(L,\cdot )$.$\spadesuit$

\begin{mycor}\label{1:9}
In a C-loop $(L,\cdot )$, if $(I,L_z^2,JL_z^2J)\in AUT(L),$ and
$(R_z^2,I,JR_z^2J)\in AUT(L)$ for all $z\in L$, then the following
are true :
\begin{enumerate}
\item $L$ is flexible. \item $(xy)^2=(yx)^2$ for all $x,y\in L$.
\item $x\mapsto x^3$ is an anti-automorphism.
\end{enumerate}
\end{mycor}
{\bf Proof}\\
This follows by Theorem~\ref{c:centralsquare}, Lemma~5.1 and
Corollary~5.2 of \cite{phd59}.$\spadesuit$

\begin{myth}\label{1:10}
A central square C-loop of exponent 4 is a group.
\end{myth}
{\bf Proof}\\
To prove this, it shall be shown that the right inner mapping

$R(x,y)=I$ for all $x,y\in L$. Corollary~\ref{1:9} is used. Let
$w\in L$.

$wR(x,y)=wR_xR_yR_{xy}^{-1}=(wx)y\cdot (xy)^{-1}=(wx)(x^2yx^2)\cdot
(xy)^{-1}=(wx^3)(yx^2)\cdot (xy)^{-1}=(w^2(w^3x^3))(yx^2)\cdot
(xy)^{-1}=(w^2(xw)^3)(yx^2)\cdot (xy)^{-1}=w^2(xw)^3\cdot
(yx^2)(xy)^{-1}=w^2(xw)^3\cdot [y\cdot x^2(xy)^{-1}]=w^2(xw)^3\cdot
[y\cdot x^2(y^{-1}x^{-1})]=w^2(xw)^3\cdot [y(y^{-1}x^{-1}\cdot
x^2)]=w^2(xw)^3\cdot [y(y^{-1}x)]=w^2(xw)^3\cdot x=w^2(w^3x^3)\cdot
x=w^2\cdot (w^3x^3)x=w^2\cdot (w^3x^{-1})x=w^2w^3=w^5=w$ if and only
if $R(x,y)=I$ if and only if $R_xR_yR_{xy}^{-1}=I$ if and only if
$R_xR_y=R_{xy}$ if and only if $zR_xR_y=zR_{xy}$ if and only if
$zx\cdot y=z\cdot xy$ if and only if $L$ is a group. Hence the claim
is true.$\spadesuit$

\begin{mycor}\label{1:11}
In a C-loop $(L,\cdot )$, if $(I,L_z^2,JL_z^2J)\in AUT(L),$ and
$(R_z^2,I,JR_z^2J)\in AUT(L)$ for all $z\in L$, then $L$ is a group.
\end{mycor}
{\bf Proof}\\
This follows from Theorem~\ref{c:centralsquare} and
Theorem~\ref{1:10}.$\spadesuit$

\begin{myrem}
Central square C-loops of exponent 4 are A-loops.
\end{myrem}

\section{Acknowledgement} The second author would
like to express his profound gratitude to the Swedish
International Development Cooperation Agency (SIDA) for the
support for this research under the framework of the Associateship
Scheme of the Abdus Salam International Centre for
theoretical Physics, Trieste, Italy. ~\\

~\\
~\\
~\\
~\\
\begin{enumerate}
\item Department of Mathematics,\\
Obafemi Awolowo University,\\ Il\'e If\`e, Nigeria.\\
{\bf e-mail}: jaiyeolatemitope@yahoo.com \item Department of Mathematics,\\
University of Ab\d e\'ok\`uta,\\ Ab\d e\'ok\`uta 110101, Nigeria.\\
{\bf e-mail}: ekenedilichineke@yahoo.com
\end{enumerate}

\begin{thebibliography}{99}
\bibitem{phd64} J. O. Ad\'en\'iran, {\it The study of properties
of certain class of loops via their Bryant-Schneider groups},
Ph.D. thesis, University of Agriculture, Abeokuta, Nigeria, 2002.
\bibitem{phd55} J. O. Ad\'en\'iran and A. R. T. Solarin, {\it A
note on generalized Bol identity}, Scientific Annals of Al.I.Cuza.
Univ., {\bf 45} 1999, 99--102.
\bibitem{phd169} A. Beg, {\it A theorem on C-loops},
Kyungpook Math. J. {\bf 17(1)} 1977, 91--94.
\bibitem{phd170} A. Beg, {\it On LC-, RC-, and C-loops}, Kyungpook Math. J. {\bf 20(2)} 1980, 211--215.
\bibitem{phd41} R. H. Bruck, {\it A survey of binary systems}, Springer-Verlag, Berlin-G\"ottingen-Heidelberg,
1966.
\bibitem{phd62} R. Capodaglio Di Cocco, {\it On Isotopism and Pseudo-Automorphism of the loops}, Bollettino U. M. I., {\bf 7} 1993, 199--205.
\bibitem{phd54} O. Chein, {\it A short note on supernuclear
(central) elements of inverse property loops}, Arch. Math., {\bf
33} 1979, 131--132.
\bibitem{phd39} O. Chein, H. O. Pflugfelder and J. D. H. Smith, {\it Quasigroups and Loops : Theory and Applications}, Heldermann Verlag,
1990.
\bibitem{phd49} J. Dene and A. D. Keedwell, {\it Latin squares and their applications}, the English University press Lts,
1974.
\bibitem{phd50} F. Fenyves, {\it Extra Loops I}, Publ. Math. Debrecen, {\bf 15} 1968, 235--238.
\bibitem{phd56} F. Fenyves, {\it Extra Loops II}, Publ. Math. Debrecen, {\bf 16} 1969, 187--192.
\bibitem{phd42} E. G. Goodaire, E. Jespers and C. P. Milies, {\it Alternative  Loop Rings}, NHMS(184), Elsevier,
1996.
\bibitem{phd60} T. G. Jaiy\'e\d ol\'a, {\it An isotopic study of
properties of central loops}, M.Sc. dissertation, University of
Agriculture, Abeokuta, Nigeria, 2005.
\bibitem{phd22} M. K. Kinyon, K. Kunen, J. D. Phillips, {\it A generalization of Moufang and Steiner loops}, Alg. Univer., {\bf 48(1)} 2002, 81--101.
\bibitem{phd59} M. K. Kinyon, J. D. Phillips and P. Vojt\v echovsk\'y , {\it C-loops : Extensions and construction}, J. Alg. \& its Appl. (to appear).
\bibitem{phd3} H. O. Pflugfelder, {\it Quasigroups and Loops : Introduction}, Sigma series in Pure Math. 7, Heldermann Verlag, Berlin,
1990.
\bibitem{phd9} J. D. Phillips and P. Vojt\v echovsk\'y,  {\it The varieties of loops of Bol-Moufang type}, Alg. Univ., {\bf 53(3)} 2005, 115-137.
\bibitem{phd61} J. D. Phillips and P. Vojt\v echovsk\'y, {\it The varieties of quasigroups of Bol-Moufang type : An equational
reasoning approach} J. Alg., {\bf 293} 2005, 17-33
\bibitem{phd58} J. D. Phillips and P. Vojt\v echovsk\'y, {\it On C-loops }, Publ. Math. Debrecen, {\bf 68(1-2)} 2006, 115-137.
\bibitem{phd10} V. S. Ramamurthi and A. R. T. Solarin, {\it On finite right central
loops}, Publ. Math. Debrecen, {\bf 35} 1988, 261--264.
\bibitem{phd53} A. R. T. Solarin, {\it On the identities of
Bol-Moufang type}, Koungpook Math. J., {\bf 28(1)} 1998, 51--62.
\bibitem{phd51} A. R. T. Solarin, {\it On certain Aktivis
algebra}, Italian Journal of Pure And Applied Mathematics, {\bf 1}
1997, 85--90.
\bibitem{phd52} A. R. T. Solarin and V. O. Chiboka, {\it A
note on G-loops}, Collections of Scientific Papers of the Faculty
of Science Krag., {\bf 17} 1995, 17--26.
\end{thebibliography}
\end{document}